\documentclass[12pt]{amsart} 

\usepackage[utf8]{inputenc} 

\usepackage{xy, graphicx, color,hyperref,verbatim,cleveref,tikz-cd}
\pdfoutput=1

\xyoption{all} 

\usepackage{amsfonts,array}

\usepackage{amssymb,amsmath,amsthm,bm,mathrsfs,enumerate,bbm,enumitem,abraces}
\usepackage[a4paper, total={6in, 8in}]{geometry}
\usepackage{fancyhdr}
\usepackage{mathtools,cite,mathdots}
\usepackage{tikz}
\usepackage{tikz-cd}
\usetikzlibrary{matrix}
\usepackage{siunitx}
\usepackage{booktabs}
\usepackage{blindtext}
\usepackage{booktabs} 

\usetikzlibrary{shapes.geometric}

\newcommand{\nc}{\newcommand}

\newcolumntype{P}[1]{>{\centering\arraybackslash}p{#1}}
\PassOptionsToPackage{linktocpage}{hyperref}
   \hypersetup{
           breaklinks=true,  
           colorlinks=true,  
           pdfusetitle=true, 
           citecolor={blue!80!black}
        }
\nc{\mc}{\mathcal}
\nc{\mb}{\mathbb}
\nc{\mf}{\mathfrak}
\nc{\ul}{\underline}
\nc{\ol}{\overline}
\nc{\dmo}{\DeclareMathOperator}
\nc{\R}{\mb R}
\dmo{\Spin}{Spin}
\dmo{\SO}{SO}
 \dmo{\pr}{pr}
 \dmo{\Sym}{Sym}
\dmo{\U}{U}
\dmo{\Hom}{Hom}
\dmo{\PGL}{PGL}
\dmo{\PSL}{PSL}
\dmo{\ortho}{orth}
\dmo{\sgn}{sgn}
\dmo{\dyn}{dyn}
\dmo{\trace}{Trace}
\dmo{\new}{new}
\dmo{\Ad}{Ad}
\dmo{\sym}{sym}
\dmo{\pal}{pal}
\dmo{\sd}{sd}
\dmo{\Ch}{Ch}
\dmo{\Span}{Span}
\dmo{\ord}{ord}
\dmo{\Or}{O}
\dmo{\Irr}{Irr}
 
\dmo{\Pin}{Pin}
\dmo{\Od}{Od}
\dmo{\EG}{EG}
\dmo{\BG}{BG}
\dmo{\ESO}{ESO}
\dmo{\BSO}{BSO}
\dmo{\BH}{BH}
\dmo{\EH}{EH}
\dmo{\Sgn}{Sgn}
\dmo{\irr}{Irr}
\dmo{\orb}{Orb}
\dmo{\rep}{Rep}
\dmo{\orep}{ORep}
\dmo{\odd}{odd}
\dmo{\dett}{det}
\dmo{\St}{St}
\dmo{\diag}{diag}
\dmo{\Top}{top}

\newcommand{\half}{\frac{1}{2}}
\newcommand{\as}{\alpha}

\newcommand{\f}{\mathbb{F}}
\newcommand{\z}{\mathbb{Z}}

\newcommand{\cc}{\mathbb{C}}

\nc{\vt}{\vartheta}

\usetikzlibrary{calc}
\tikzset{
between/.style args={#1 and #2}{
    at = ($(#1)!0.5!(#2)$)
},
betweenl/.style args={#1 and #2}{
    at = ($(#1)!0.35!(#2)$)
}
}

\newenvironment{roster}
 {\begin{enumerate}[font=\upshape,label=(\roman*)]}
 {\end{enumerate}}

 \newtheorem{thm}{Theorem}[section]
\newtheorem{c.intro}[thm]{Corollary}
\newtheorem{lemma}[thm]{Lemma}
\newtheorem{prop}[thm]{Proposition}

\newtheorem{cor}{Corollary}[thm]

\theoremstyle{definition}

\newtheorem{ex}[thm]{Example}
\theoremstyle{definition}
\newtheorem{remark}[thm]{Remark}
\theoremstyle{definition}

\dmo{\Mod}{mod}
\dmo{\res}{res}
 \dmo{\Sq}{Sq}
  \dmo{\Tr}{Tr}
 \dmo{\RO}{RO}
\dmo{\Sp}{Sp}
\dmo{\SL}{SL}
\dmo{\GL}{GL}
\dmo{\GSp}{GSp}
\nc{\la}{\lambda}
\nc{\eps}{\varepsilon}
\nc{\lip}{\langle}
 \nc{\rip}{\rangle}
\nc{\gm}{\gamma}
\nc{\beq}{\begin{equation*}}
\nc{\eeq}{\end{equation*}}
\dmo{\Perm}{Perm}
\dmo{\Res}{Res}
\dmo{\Ind}{Ind}
\dmo{\ind}{ind}
\dmo{\tr}{tr}
\dmo{\reg}{reg}
\dmo{\End}{End}
\dmo{\SW}{SW}
\dmo{\Syl}{Syl}
\MakeRobust\(
\MakeRobust\)

\title[Stiefel-Whitney Classes]{Stiefel-Whitney Classes for finite special linear groups of even rank }
\author{Neha Malik}
\author{Steven Spallone}
\date{}
\address{Chennai Mathematical Institute, India}
\email{51nehamalik94@gmail.com}
\address{Indian Institute of Science Education and Research, Pune-411008,Maharashtra,India}
\email{sspallone@gmail.com}
\begin{document}
\maketitle

\begin{abstract}
We compute the total Stiefel-Whitney Classes (SWCs) for orthogonal representations of special linear groups $\SL(n,q)$ when $n$ and $q$ are odd. These classes are expressed in terms of character values at diagonal elements of order $2$. We give several consequences, and  work out the 4th SWC explicitly, and the 8th SWC when the 4th vanishes.
\end{abstract}
\tableofcontents
\section{Introduction} 
Stiefel-Whitney Classes (SWCs) are important cohomological invariants of (orthogonal complex) vector bundles. This paper is part of a study of SWCs of vector bundles arising from orthogonal representations of well-known groups. SWCs for $\GL(n,q)$ were treated in \cite{GJgln}, and in \cite{NSSL2} we studied $\SL(2,q)$. Earlier work on computing SWCs includes \cite{guillot}, which in particular presents a computer program that treats $2$-groups up to order $64$. Appendix A of the same paper surveys  earlier work in this topic.

Let $q$ be an odd prime power, $n\geq 3$ be odd, and $G=\SL(n,q)$. Let $\pi$ be an orthogonal representation of $G$. We show how to compute  each SWC $w_k(\pi)$ in principle, and perform this calculation explicitly for some small $k$.
 We express the total SWC $w(\pi)$ as a product, from which each $w_k(\pi)$ can be calculated in principle.  
 
We factor $w(\pi)$ into a product of certain symmetric polynomials;  the exponents of this factorization are written in terms of character values of $\pi$ on diagonal elements of order $2$.  To find an individual SWC $w_k(\pi)$, one only needs these polynomials up to degree $k/2$ or $k$, depending on whether $q$ is $1$ or $3 \pmod 4$.

The first three (nontrivial) SWCs of $\pi$ automatically vanish. Generally, the degree of the first nonvanishing SWC is a power of $2$. We explicitly compute $w_4(\pi)$, and when $q \equiv 1 \pmod 4$, we also give $w_8(\pi)$. 
To mention here a resulting vanishing condition, put   $b_2=\diag(-1,-1,1,\hdots, 1)$ and $b_4=\diag(-1,-1,-1,-1,1,\hdots,1)\in G$. 
Let  \begin{align*}
m_1&=\frac{1}{16}\big(\chi_\pi(\mathbbm{1})+2\chi_\pi(b_2)-3\chi_\pi(b_4)\big) \text{ and }\\
m_2&=\frac{1}{16}\big(\chi_\pi(\mathbbm{1})-2\chi_\pi(b_2)+\chi_\pi(b_4)\big).
\end{align*}
In fact, $m_1 \equiv m_2 \pmod{2}$. When $q \equiv 1 \pmod 4$, these are both even.

\begin{thm} \label{intro.thm}
  Let $\pi$ be an orthogonal representation of $\SL(n,q)$ with odd $n \geq 5$. Then $w_4(\pi)=0$ iff $m_1 \equiv m_2 \pmod 4$. Moreover, when   $q \equiv 1 \pmod 4$ and  $w_4(\pi)=0$, then  $w_8(\pi)=0$ iff $m_1 \equiv m_2 \pmod 8$.

\end{thm}
 
 In fact, it is enough to do the $w_4$-calculation for $G=\SL(5,q)$ when $q \equiv 3 \pmod 4$, since the restriction map on (degree $4$) cohomology is injective.  To find a universal formula for $w_8(\pi)$ when $q \equiv 3 \pmod 4$, or for $w_{16}(\pi)$ when $q \equiv 1 \pmod 4$, one would need to compute these  in $\SL(9,q)$. Similarly, one can make ``universal'' calculations for any $w_k$ with $k$ even, by working with $\SL(k+1,q)$. (See Remark \ref{univ.remark}.)

For $G=\SL(3,q)$ we also determine the subring $H^*_{\SW}(G,\z/2\z)$ of $H^*(G,\z/2\z)$ generated by SWCs $w_k(\pi)$  , the subgroup $\mathcal W(G)$ generated by the total SWCs $w(\pi)$, and for each $\pi$ determine the smallest nontrivial $w_k(\pi)$.

The layout of the paper is as follows: Section \ref{not} sets up notations and recalls some fundamental notions used throughout the paper. In some sense, our calculation reduces to a study of $S_n$ and $S_{n+1}$-invariant representations of elementary abelian $2$-groups, which we consider in Section \ref{eab2}. 
In Section \ref{revGJ}, we review the calculation of SWCs for general linear groups in \cite{GJgln}. In Section \ref{slg}, we compute the total SWCs for $\SL(n,q)$.  For $G=\SL(3,q)$, the calculation of obstruction classes, $H^*_{\SW}(G,\z/2\z)$ and the top SWC are also carried out in Section \ref{sl33}. We also treat the group $\SL(5,q)$ in Section \ref{sl5}. In the final Section \ref{uniswc}, we illustrate the method of universal calculations, by computing $w_4$ and $w_8$ for all $\SL(n,q)$.

\textbf{Acknowledgments.} 
The authors would like to thank Rohit Joshi for helpful discussions, and Varun Shah for the calculations of Section \ref{dsonfac}. Part of this paper comes out of the first author’s Ph.D. thesis \cite{Malik.thesis} at IISER Pune, during which she was supported by a Ph.D. fellowship from the Council of Scientific and Industrial Research, India. Both authors thank the Chennai Mathematical Institute for its hospitality during many visits.

\section{Notations and preliminaries}\label{not}

\subsection{Representations}

Let $G$ be a finite group. We will only consider finite dimensional complex representations $(\pi,V)$ of $G$.
  Write $(\pi^\vee,V^\vee)$ for the dual representation.  Write $\chi_\pi(g)$ for the character value of $\pi$ at $g\in G$. A \emph{linear character} $\chi$ is a degree $1$ representation, and  $\chi$   is \emph{quadratic} when $\chi(g)=\pm 1$ for all $g$. If $\pi_1$ is a representation of $G_1$ and $\pi_2$ is a representation of $G_2$, write $\pi_1 \boxtimes \pi_2$ for the external tensor product representation of $G_1 \times G_2$.

When $V$ admits a $G$-invariant symmetric nondegenerate bilinear form, we say that $\pi$ is \emph{orthogonal}. A linear character $\chi$ is orthogonal if and only if it is quadratic. 
Given any representation $\pi$, one can form the orthogonal representation $S(\pi)=\pi \oplus \pi^\vee$ on $V \oplus V^\vee$, by defining the inner product $( (v,v^*),(w,w^*))=\lip v^*,w\rip+\lip w^*,v \rip$.  An orthogonal representation $\pi$ is \emph{orthogonally irreducible} (an OIR), when it cannot be decomposed as the sum of two orthogonal representations. There are two kinds of OIRs: irreducible representations which are orthogonal, and representations of the form $S(\pi)$ when $\pi$ is irreducible but not orthogonal. Every orthogonal representation is decomposable into the sum of OIRs.

\subsection{Stiefel-Whitney Classes}
In this section we briefly review notations for SWCs of representations; please see \cite[Section 2.3]{NSSL2} for a more generous introduction. Throughout this paper write $H^*(G)$ for the group cohomology ring $H^*(G,\z/2\z)$.

Let $\pi$ be an orthogonal representation of   $G$ of degree $d$. Associated to $\pi$ are cohomological invariants
$$w_k(\pi)\in H^k(G) \: \: \text{for} \: \:  k=0,1,2,\hdots, d$$
known as the $k$th \emph{Stiefel-Whitney Class} (SWC) of $\pi$. Their sum 
\beq
w(\pi)=w_0(\pi)+w_1(\pi)+\hdots
\eeq
 is called the \emph{total SWC} of $\pi$. One takes $w_0(\pi)=1$ and $w_k(\pi)=0$ for $k>d$. When $d=1$, 
the first SWC $w_1$ comes from the isomorphism $\Hom(G, \{ \pm 1\}) \cong H^1(G)$. More generally,
$w_1(\pi):=w_1(\det \pi)$. When $\det \pi$ is trivial, the second SWC $w_2(\pi)$ vanishes iff the representation $\pi: G \to \SO(V)$ lifts to the spin group.

  The image of the Euler class $e(\pi) \in H^{d}(G,\z)$   under the  coefficient map to $H^d(G)$ is $w_d(\pi)$, by \cite[Property 9.5]{milnor}. We will denote this class by $w_{\Top}(\pi)$. 
Let $H^\bullet(G)=\prod_i H^i(G)$ be the complete mod $2$ cohomology ring of $G$, consisting of all formal infinite series $\alpha_0+\alpha_1+ \cdots$, with each $\alpha_i \in H^i(G)$. A series with $\alpha_0=1$ is invertible in $H^\bullet(G)$, and therefore we may define $w(\pi)$ for virtual representations of $G$ by $w(\pi \ominus \pi')=w(\pi) \cup w(\pi')^{-1} \in H^\bullet(G)$, when $\pi$ and $\pi'$ are orthogonal representations.

Write $\mc W(G)$ for the ``Stiefel-Whitney group", meaning the multiplicative subgroup of $H^\bullet(G)$ generated by the total SWCs $w(\pi)$ of orthogonal representations. Write $H_{\SW}^*(G)$ for the subring of $H^*(G)$ generated by SWCs $w_k(\pi)$ of orthogonal representations $\pi$. 
 
 \begin{lemma} \label{Schur} Let $G=\SL(n,q)$ for $n\geq 3$ a positive integer and $q$ odd. For any orthogonal representation $\pi$ of $G$, 
 we have $w_1(\pi)=w_2(\pi)=0$.
 \end{lemma}
 
 \begin{proof} Since $G$ is perfect, $\det \pi$ is trivial. But $G$ also has trivial Schur multiplier and therefore $\pi$ lifts to the spin group, e.g.,  by \cite[Proposition 6]{joshi2021spinoriality}.  
 \end{proof}
 
 Let $\pi$ be an orthogonal representation of $G$, and $G' \leq G$  a subgroup. Write $w^{G'}(\pi)$ be the SWC of the restriction of $\pi$ to $G'$.
 
\subsection{Detection} 

An important technique in this subject is cohomological detection. A subgroup $G' \leq G$ \emph{detects the mod $2$ cohomology of $G$}, when the restriction map $i^*: H^*(G) \to H^*(G')$ is injective. Moreover we say $G'$ \emph{detects SWCs of $G$}, when the restriction of $i^*$ to $H^*_{\SW}(G)$ is injective. 
For example, a Sylow 2-subgroup of $G$ detects its mod 2 cohomology, e.g., by  \cite[ Cor. 5.2, Ch. II]{Milgram}.

\begin{thm}[\cite{Quillen}, Theorem 3]\label{quill}
 The diagonal subgroup $T$ detects the mod $2$ cohomology of $\GL(n,q)$.
\end{thm} 

Let $N=N_G(G')$ be the normalizer in $G$ of $G'$, and put $W=N/G'$. Then the restriction of $\pi$ to $G'$ is a $W$-invariant representation, and the image of $i^*$ lies in the fixed points $H^*(G')^W$. For instance, the image of $H^*(\GL(n,q))$ lies in $H^*(T)^{W}$, where $W$ is the subgroup of permutation matrices in $\GL(n,q)$. 
 
 When $G'$ detects the mod $2$ cohomology or SWCs of $G$, and $\pi$ is an orthogonal representation of $G$,  we may write $w(\pi)$ for $w^{G'}(\pi)$ by abuse of notation.

\subsection{Symmetric Functions}
 
The elementary symmetric functions $\mc{E}_k(\bm x)$ over $\z/2\z$ in variables $x_1,x_2,\hdots,x_n$ are defined for $k\leq n$ as,
$$\mc{E}_k(\bm x):=\sum_{1\leq i_1<\hdots<i_k\leq n}x_{i_1}x_{i_2}\hdots x_{i_k}\in \z/2\z[x_1,\hdots,x_n].$$
(By convention these are $0$ for $k>n$.)
For example, $\mc E_1(x_1,x_2)=x_1+x_2$ and $\mc E_2(x_1,x_2)=x_1x_2$. We have
$$\prod_{i=1}^n (1+x_i)=1+\mc E_1(\bm x)+\hdots+\mc E_n(\bm x).$$
Similarly, there are complete symmetric functions $\mc H_k(\bm x) \in \z/2\z[x_1,\hdots,x_n]$ which are defined as,
$$\mc{H}_k(\bm x):=\sum_{1\leq i_1\leq\hdots\leq i_k\leq n}x_{i_1}x_{i_2}\hdots x_{i_k}\quad \text{for}\quad k\leq n$$
and $0$ for $k> n$. For example, $\mc H_1(x_1,x_2)=x_1+x_2$ and $\mc H_2(x_1,x_2)=x_1^2+x_2^2+x_1x_2$.

\subsection{Cyclic Groups}\label{cy}
Let $C=C_\ell$, the cyclic group of even order $\ell$, and fix a generator $g$. We set $-1=g^{\ell/2}$, the unique order $2$ element of $C$. Write $\sgn:C \to \{ \pm 1\}$ for the linear character with $\sgn(g)=-1$.
A linear character $\psi$  of $C$ is \textit{even} when $\psi(-1)=1$, and  \textit{odd} when $\psi(-1)=-1$.   
Let $\psi_\bullet$ be the linear character  of $C$ with $\psi_\bullet(g)=e^{\frac{2\pi i}{\ell} }$. It is known \cite{KT} that 
$$H^*(C)=\begin{cases}
\z/2\z[s,t]/(s^2),& \ell\equiv0 \text{ (mod } 4)\\
\z/2\z[v],& \ell\equiv2 \text{ (mod }4)
\end{cases}$$where $s=w_1(\sgn),$ $t=w_2(S(\psi_\bullet))$ for $\ell\equiv0 \text{ (mod } 4)$, and $v=w_1(\sgn)$ when $\ell \equiv2 \text{ (mod }4)$.

\begin{prop}\label{cyswc} Let $\pi$ be an orthogonal representation of $C$. Put $b_\pi=\frac{1}{2} \left( \deg \pi-\chi_\pi(-1) \right)$.
\begin{roster}
\item If  $\ell\equiv 2\: (\Mod\: 4)$, then $w(\pi)=(1+v)^{b_\pi}$.
\item If  $\ell\equiv 0\: (\Mod\: 4)$, then $b_\pi$ is even and $w(\pi)=(1+\delta_\pi s)(1+t)^{b_\pi/2}$,
 where  
 $$\delta_\pi=\begin{cases}
0,&\text{ if } \det \pi=1\\
1,& \text{ if } \det \pi=- 1.
\end{cases}$$
\end{roster}
\end{prop}
 
 \begin{proof} Straightforward; see   \cite[Proposition 2.14]{Malik.thesis}. \end{proof}

Let $C^r$ be the $r$-fold product of $C$, with projection maps $\pr_i: C^r\to C$ for $1 \leq i \leq r$. By K\"{u}nneth, we have 
\begin{equation}\label{Cnk}H^*(C^r)=\begin{cases}
\z/2\z[s_1,\dots,s_r,t_1,\dots,t_r]/(s_1^2,\hdots,s_r^2),& \ell\equiv0 \text{ (mod } 4)\\
\z/2\z[v_1,\dots,v_r],& \ell\equiv2 \text{ (mod }4)
\end{cases}\end{equation}
where we put $s_i=w_1(\sgn\circ\pr_i)$ and $t_i=w_2(S(\psi_\bullet)\circ\pr_i)$  for $\ell\equiv0\:(\Mod\:4)$, and $v_i=w_1(\sgn\circ\pr_i)$, for $\ell\equiv2 \text{ (mod } 4)$.

We will later need the following:
\begin{lemma} \label{about.P} Let $\ell$ be a multiple of $4$. The restriction map $H^2(C^r) \to H^2(\{ \pm 1\}^r)$ takes $t_i$ to $v_i^2$.
\end{lemma} 

\begin{proof} This amounts to the fact that the restriction of $\psi_\bullet$ to $\{ \pm 1\}$ is nontrivial.
\end{proof}

\begin{lemma} \label{lien} Let $\ell$ be a multiple of $4$, and $\pi$ an orthogonal representation of $A=C_\ell^r$.
The multiplicity of a nontrivial linear character of $E=C_2^r$ in the restriction of $\pi$ to $E<A$ is even.
\end{lemma}

\begin{proof} We may assume $\pi$ is orthogonally irreducible. It is easy to see that an orthogonal linear character of $A$ restricts to the trivial representation of $E$. Note that if $\psi$ is a representation of $A$, then $\Res^A_E \psi \cong \Res^A_E \psi^\vee$. It follows that the multiplicity of a nontrivial linear character of $E$ in $\Res^A_E S(\psi)$ is even.
\end{proof} 

\section{Elementary Abelian $2$-groups}\label{eab2}

\subsection{Linear Characters}

Let $E=C_2^r$,  an elementary abelian $2$-group. All irreducible representations of $E$ are quadratic linear characters. Since $H^*(E)$ is a polynomial ring \eqref{Cnk}, the largest $k$ with $w_k(\pi) \neq 0$  is the number of nontrivial irreducible constituents of $\pi$.
In particular:
\begin{prop} \label{E.topSWC} The class $w_{\Top}(\pi) \neq 0$ iff $\pi$ does not contain the trivial representation.
\end{prop}

View $E$ as an $\mb F_2$-vector space, say with basis $e_1,\hdots,e_r$. For $e\in E,$ we write $|e|=\#\{i:c_i=1\}$ when $e=\sum\limits_{i=1}^r c_ie_i$ with $c_{i}\in \mathbb F_2$. Let $E^\vee=\Hom(E,\mathbb{F}_2)$ and consider the basis $v_1, \hdots, v_r$ of $E^\vee$ dual to the $e_i$.  For $v\in E^\vee,$ put $|v|=\#\{i:c_i=1\}$ when $v=\sum\limits_{i=1}^r c_iv_i$ with $c_{i}\in \mathbb F_2$.

\subsection{$S_r$-invariant Representations of $E$}\label{eab}
The symmetric group $S_r$ acts on $E$ by permuting the $e_i$. For $0 \leq k \leq r$, put
\beq
\mathcal O_k=\Big\{e\in E: |e|=k\Big\};
\eeq
these are the $S_r$-orbits in $E$. For $1 \leq k \leq r$, a representative is given by $a_k=e_1+ \cdots + e_k  \in \mc O_k$, and we put $a_0=0$.
 Identifying $E^\vee$ with $H^1(E)$, we may write
$$H^*(E)=\Sym(E^\vee)\cong \z/2\z[v_1,\hdots,v_r]$$
where $\Sym(E^\vee)$ is the symmetric algebra of $E^\vee$.  The dual of above action of $S_r$ on $E^\vee$ permutes the dual basis. Now the $S_r$-orbits in $E^\vee$ are the sets
$$\mathcal{O}^*_k 
=\{v\in E^\vee:|v|=k\},$$
for $0 \leq k \leq r$. A representative is given by $\vt_k=v_1+v_2+\hdots + v_k \in \mc O_k^*$, and $\vt_0=0$.
Note that $\{ \vt_1, \ldots, \vt_r \}$ is a basis for $E^\vee$.

The representation $$\sigma_k=\bigoplus_{v \in \mc O_k^*} v \quad;\quad k=0,1,\hdots,r$$
 is $S_r$-invariant, of degree $\binom{r}{k}$. (Here we view each $v \in \mc O_k^*$ as a linear character $E \to \cc^\times$, and so $\sigma_k$ is a complex representation.) For instance, $\sigma_0$ is the trivial linear character.
Any $S_r$-invariant representation $\sigma$ of $E$  decomposes uniquely into a direct sum of $\sigma_k$'s, so we may speak of the multiplicity $m_k(\sigma)$ of $\sigma_k$ in $\sigma$. (Note this is also the multiplicity of $\vt_k$ in $\sigma$.) In other words,
\begin{equation}\label{pisig}\sigma=\bigoplus_{k=0}^r m_{k}(\sigma)\sigma_k.\end{equation}
These multiplicities were computed in  \cite{GJgln}, in terms of character values of $\sigma$.  Then
\begin{equation}\label{mk}m_k(\sigma)=\frac{1}{2^{r}}\sum\limits_{i=0}^r\chi_{\sigma_i}(a_k)\chi_\sigma(a_i),\end{equation}
and we are left with determining the constants $\chi_{\sigma_i}(a_k)$.
 For a polynomial $f$, write $[f]_i$ for the coefficient of the degree $i$ term of $f$; in other words so that $f(x)=\sum_i [f]_i x^i$.

\begin{lemma}[\cite{GJgln}, Proposition 2]\label{chisig}
We have 
\beq
\chi_{\sigma_i}(a_k)= \left[ (1-x)^k (1+x)^{r-k} \right]_i.
\eeq
\end{lemma}
For example:
\begin{equation} \label{m0.formula}
m_0(\sigma)=\frac{1}{2^r} \sum_{i=0}^r \binom{r}{i} \chi_\sigma(a_i).
\end{equation}
 
To summarize: When $E$ is an elementary abelian $2$-group of rank $r$ with a given (unordered) basis, then one has the fundamental symmetric representations $\sigma_k$ as above, and the multiplicity of $\sigma_k$ in an $S_r$-invariant representation $\sigma$ is calculated by \eqref{mk} and Lemma \ref{chisig}.
 
\subsection{$S_{r+1}$-invariant representations of $E$}\label{Winv}
 Now suppose that the rank of $E$ is even, say $r=2\ell$. We can identify $E$ with the hyperplane
 $E^+_0=\{(a_1,\hdots,a_{r+1}): \sum_{i=1}^{r+1} a_i=0\}$ in $E^+=\mathbb{F}_2^{r+1}$ via $\sum a_i e_i \mapsto (a_1,a_2,\hdots,a_r,a_1+\hdots+a_r)$. Since $E^+_0$ is fixed by $S_{r+1}$, this gives an action of $S_{r+1}$ on $E$ which further induces an action of $S_{r+1}$ on $E^\vee$. 
\begin{lemma} We have,
\begin{roster}
\item The orbits of $E$ under the action of $S_{r+1}$:
$$\mc O_0, \mc O_1 \cup \mc O_2, \hdots, \mc O_{2i-1}\cup \mc O_{2i},\hdots,  \mc O_{r-1}\cup \mc O_{r}. $$
\item The orbits of $E^\vee$ under $S_{r+1}$:
$$
\mathcal{O}_0^*, \mathcal{O}_1^*\cup \mathcal{O}_{r}^*,\hdots, 
\mathcal{O}_k^*\cup \mathcal{O}_{r+1-k}^*,\hdots, \mathcal{O}_\ell^*\cup \mathcal{O}_{r+1-\ell}^*.
$$
\end{roster}
\end{lemma}

\begin{proof}
($1$) is elementary. Write $\vt \sim \vt'$ when $\vt,\vt'$ are in the same $S_{r+1}$-orbit. Since $\vt_k \in \mc O_k^*$, for ($2$) it is enough to show that 
 $\vt_i \sim \vt_j$ are in the same $S_{r+1}$-orbit iff $j=i$ or $j=r+1-i$ for non-zero $i,j$.
 Let $\tau=(1,r+1)\in S_{r+1}.$ One can check 
$$\tau\cdot\vt_k=v_1+v_{k+1}+\hdots+v_r \in \mc O^*_{r+1-k}$$
so that $\vt_k \sim \vt_{r+1-k}$. Now by a result of Brauer (see \cite[Cor. 6.33]{isch}), the number of $S_{r+1}$-orbits 
in $E^\vee$ is the same as the number of $S_{r+1}$-orbits in $E$, namely $\ell+1$. This forces  $(2)$.
\end{proof}

\begin{cor}\label{minrep} The representations $\widetilde \sigma_k=\sigma_k \oplus \sigma_{r+1-k}$ with $1 \leq l \leq \ell$, together with $\sigma_0$, are the minimal $S_{r+1}$-invariant representations of $E$. \end{cor}

\subsection{Dickson Factors}\label{dsonfac}
We define polynomials $\mc D^{[k]}=\mc D_r^{[k]}\in \mathbb{F}_2[x_1,\hdots,x_r]$ as,
\beq
 \mc D^{[k]}(\bm x)=\prod_{1\leq i_1<\hdots<i_k\leq r}(1+x_{i_1}+\hdots+x_{i_k}),
\eeq
where $\bm x$ is the unordered tuple $\{x_1, \ldots, x_r\}$. Also put $\mc D(\bm x):=\prod_{k} \mc D^{[k]}({\bm x})$. For example, if $r=3$, then
\beq
\mc D^{[2]}({\bm x})= (1+x_1+x_2)(1+x_1+x_3)(1+x_2+x_3). \eeq
Let $E$ be as in the previous section, and put ${\bm v}=\{v_1, \ldots, v_r\}$. In this paper, we will encounter expressions of the form $\prod\limits_{|v|=k}(1+v)\in \Sym(E^\vee)$. We note that 
\beq
\prod_{|v|=k }(1+v)= \mc D^{[k]}({\bm v}).
\eeq

The homogeneous terms of  $\mc D(\bm x)$ are well-understood: we have  $\mc D(\bm x)=1+\sum_{i=1}^r d_i(\bm x)$,
where $d_i$ are certain $\GL(r,2)$-invariant polynomials of  degree $2^r-2^{r-i}$,  known as \emph{Dickson invariants}. For elaboration, see \cite[ Section 5.2]{NSSL2} and \cite{wilk}. We call $\mc D^{[k]}$, the $k$th \emph{Dickson factor}. 
 
Clearly, $\mc D^{[k]}$ is a symmetric polynomial. We have $\mc D^{[1]}=1+\mc E_1+ \mc E_2+ \cdots + \mc E_r$,
where $\mc E_i=\mc E_i({\bm x})$   is the $i$th elementary symmetric polynomial in the $x_i$.
Also, $\mc D^{[r]}=1+\mc E_1$.
So for $r=2$, 
\beq
\mc D=\mc D^{[1]}\mc D^{[2]}=1+ \underbrace{\mc E_1^2+\mc E_2}_{d_1}+\underbrace{\mc E_1 \mc E_2}_{d_2}.
\eeq

If $r=3$, then
\beq
\mc D^{[2]}= 1+(\mc E_1^2+\mc E_2)+( \mc E_1 \mc E_2+ \mc E_3). \\
\eeq

If $r=4$, then
\beq
\mc D^{[2]}=1+ \mc E_1+ \mc E_1^2 + \mc E_1^3+( \mc E_2^2+\mc E_1\mc E_3)+(\mc E_1\mc E_2^2+\mc E_1^2\mc E_3)+(\mc E_1\mc E_2\mc E_3+\mc E_3^2+\mc E_1^2\mc E_4)
\eeq
and
\beq
\mc D^{[3]}=1+ \mc E_1+ (\mc E_1^2 + \mc E_2)+( \mc E_3+ \mc E_1^3)+( \mc E_1^2\mc E_2+ \mc E_1 \mc E_3+\mc E_4).
\eeq

\section{A review of Joshi-Ganguly's work on $\GL(n,q)$}\label{revGJ}
The SWCs of orthogonal representations of $G=\GL(n,q)$ were found in \cite{GJgln}.  By Theorem \ref{quill},  $H^*(G)$ injects into $H^*(T)^{S_n}$, where  $T<G$ is the diagonal subgroup. Write $T[2] <T $ for the subgroup of diagonal matrices, all of whose eigenvalues are $\pm 1$. It is an elementary abelian $2$-group of rank $n$, with basis \newline
$e_i=\diag(1,\hdots,1,\aunderbrace{-1}_{i},1,\hdots,1)$.

\subsection{Case $q\equiv 3\;(\Mod\: 4)$}
Here, $T[2]$ is the Sylow $2$-subgroup of $T$, hence detects the mod $2$ cohomology of $G$. 
Let ${\bm v}=\{v_1, \ldots, v_n\}$ be the basis of $E^\vee=H^1(E)$ dual to the $e_i$.
 Suppose $\pi$ is an orthogonal representation of $G$. Then  $\sigma: =\res^G_{T[2]}\pi$ is $S_n$-invariant. Following \cite[Theorem 1]{GJgln} we have
\begin{equation}\label{q34}
\begin{split}
w(\pi) &= 
\prod\limits_{k=1}^n\Big(\prod\limits_{|v|=k }(1+ v)\Big)^{m_k(\pi)} \\
&=  \prod\limits_{k=1}^n\mc D_n^{[k]}({\bm v})^{m_{k}(\pi)},
\end{split}
\end{equation}
 with $m_k(\pi)=m_k(\sigma)$ as in   Section \ref{eab}. Then, $w^{T[2]}(\pi)$ is given by \eqref{q34} in the case of $ q\equiv  3\: (\Mod\: 4)$, since $T[2]$ is a detecting subgroup.

\subsection{Case $q\equiv 1\;(\Mod\: 4)$}\label{glnq14}
 In this case, $T\cong C_{q-1}^n$ is the detecting subgroup. From \eqref{Cnk}, the mod $2$ group cohomology of $T$ is of the form
$H^*(T)\cong \Lambda \otimes \mathcal{P}$,
where $\Lambda$ is the exterior algebra $\z/2\z[s_1,\hdots,s_n]/(s_1^2,\hdots,s_n^2)$ and $\mathcal{P}$ is the polynomial algebra $\z/2\z[t_1,\hdots,t_n]$.   
For an orthogonal representation $\pi$ of $G$, \cite[Theorem 1]{GJgln} gives 
\begin{equation}\label{q14}
\begin{split}
w(\pi)&= \left(1+\delta \cdot \mc E_1({\bm s})\right)\prod\limits_{k=1}^n {\mc D}_n^{[k]}({\bm t})^{m_k(\pi)/2}.\\
\end{split}
\end{equation}
 Here, $\delta=0$ if $\det\pi=1$ and  otherwise $\delta=1$, also ${\bm s}=\{s_1, \ldots, s_n\}$ and  ${\bm t}=\{t_1, \ldots, t_n\}$.

\

\section{Special linear groups of even rank}\label{slg}
 Let $n\geq 3$ be odd, and $G=\SL(n,q)$. (Throughout this paper, $q$ is an odd prime power.)

\subsection{Detection}\label{SLND}
 Let
$$M=\Bigg\{\begin{pmatrix}
A&&0\\
0&& \det(A)^{-1}
\end{pmatrix}: A\in \GL(n-1,q)\Bigg\},$$
isomorphic to $\GL(n-1,q)$.
We view $T$ and $T[2]$ as subgroups of $G$ via $M$.  For $1 \leq i <n$, put
$e_i^+=\diag(1,\hdots,1,\aunderbrace{-1}_{i},1,\hdots,1,-1) \in \SL(n,q)$; we use this as a basis of the 
elementary abelian $2$-group $T[2] < \SL(n,q)$.

\begin{lemma}\label{dl}
The diagonal subgroup $T$ detects the mod $2$ cohomology of $G$. 
\end{lemma}

\begin{proof}
The subgroup
 $M$ has odd index in $G$, hence contains a Sylow $2$-subgroup.   Thus $\GL(n-1,q)$ detects the mod $2$ cohomology of $G$. This, with Theorem \ref{quill}, gives the injectivity of the composition $H^*(G)\hookrightarrow H^*(M)\hookrightarrow H^*(T)$.
\end{proof}

\begin{remark} If either $n$ or $q$ is even, then $[G:M]$ is even, so this argument breaks down in those cases.
\end{remark}


\begin{lemma}
The subgroup $T[2]$ detects the SWCs of $G$.
\end{lemma}

\begin{proof}When $q \equiv 3 \pmod 4$, the lemma is clear since $T[2]$ is the $2$-Sylow of $T$. So assume $q \equiv 1 \pmod 4$. Let $\pi$ be an orthogonal representation of $G$. Since $G$ is perfect, $\delta=0$ in formula \eqref{q14}, hence $w^T(\pi) \in \mc P$. But by Lemma \ref{about.P},
the restriction map 
$H^*(T) \to H^*(T[2])$ is injective on $\mc P$.
\end{proof}

Combining this with Proposition \ref{E.topSWC} gives:

 \begin{cor} \label{G.topSWC} Let $\pi$ be an orthogonal representation of $G$. The following are equivalent:
 \begin{roster}
 \item The trivial representation of $T[2]$ does not occur in the restriction of $\pi$ to $T[2]$.
 \item $m_0(\pi) =0$.
 \item $w_{\Top}(\pi) \neq 0$.
 \end{roster}
  \end{cor}
  (Here $m_0(\pi)$ is given by Equation \eqref{m0.formula}.)

\subsection{SWC Calculations}


Let $\pi$ be an orthogonal representation of $G$. By restriction we may regard it as an orthogonal representation of
$\GL(n-1,q)<G$, which is a detecting subgroup. Therefore, $w(\pi)$ is given by \eqref{q34}, \eqref{q14} above.
These formulas used the fact that the restriction of $\pi$ to $T$ is $S_{n-1}$-invariant.
 But for us, $\pi$ comes from the bigger group $\SL(n,q)$, and so this  restriction is invariant under the larger Weyl group $S_{n}$.
In this section, we exploit this larger symmetry to give a simpler formula for $w(\pi)$. 

Let $m_k(\pi)$ be the multiplicity of $\sigma_k$ in $\sigma=\pi|_{T[2]}$  as in Section \ref{eab} (with respect to the basis $b_i$). For $1 \leq i \leq \frac{n-1}{2}$, we define
$\ol \sigma_i=\sigma_{2i-1} \oplus \sigma_{2i}$. The diagonal matrices 
\beq b_k=\diag(\aunderbrace{-1,\hdots,-1}_{k},\underbrace{1,\hdots,1}_{n-1-k},(-1)^k) 
\eeq
corresponds to the elements $a_k$ of Section \ref{eab}.
  
 \begin{lemma} \label{mkpi} For $1 \leq k \leq \frac{n-1}{2}$, we have
 \begin{equation} \label{office}
 m_k(\pi)=m_{n-k}(\pi)
 \end{equation}
and
 \beq
 m_k(\pi)= \frac{1}{2^{n-1}} \left( \deg \pi+ \sum_{i=1}^{\frac{n-1}2} \chi_{\ol \sigma_i}(b_k) \chi_{\pi}(b_{2i}) \right).
 \eeq
Moreover, the character value $\chi_{\ol \sigma_i}(b_{k})$ is the coefficient of $x^{2i}$ in 
 \beq
 (1-x)^{k}(1+x)^{n-k}.
 \eeq
 \end{lemma}

 \begin{proof} The restriction of $\pi$ to $T[2]$ is $S_{n}$-invariant.  The minimal $S_{n}$-invariant representations are of the form $ \widetilde \sigma_k=\sigma_k \oplus \sigma_{n-k}$
from Corollary \ref{minrep}.
Therefore both sides of \eqref{office} equal the multiplicity of $\widetilde \sigma_k$ in $\pi$.

 From Section \ref{eab}, we have
 \beq
 m_k(\pi)= \frac{1}{2^{n-1}} \left(\sum_{i=0}^{n-1} \chi_{\sigma_i}(b_k) \chi_{\pi}(b_{i}) \right).
 \eeq
 Now $b_{2i}$ and $b_{2i-1}$ are conjugate in $G$ for $1 \leq i \leq \frac{n-1}{2}$, so  $\chi_{\pi}(b_{2i})=\chi_{\pi}(b_{2i-1})$, which gives  the first statement. For the last statement, note that
 \beq
 \begin{split}
 \chi_{\ol \sigma_i}(b_{k}) &= \left[ (1-x)^k(1+x)^{n-1-k} \right]_{2i-1}+  \left[ (1-x)^k(1+x)^{n-1-k} \right]_{2i} \\
 					&= \left[ (1-x)^k(1+x)^{n-k} \right]_{2i}.
 \end{split}
 \eeq \end{proof}
Similarly, for $k=0$ the formula \eqref{m0.formula} simplifies to
 \begin{equation} \label{m0.simp}
 m_0(\pi)= \frac{1}{2^{n-1}}  \sum_{i=0}^{\frac{n-1}{2}} \binom{n}{2i} \chi_{\pi}(b_{2i}).
 \end{equation}
 
  Now we come to the central result.

\begin{thm}\label{SLnSWCs}
 Let $\pi$ be an orthogonal representation of $G=\SL(n,q)$.
  \begin{roster}
 \item For  $q \equiv 1 \pmod 4$,
\beq w(\pi)=\prod\limits_{k=1}^{\frac{n-1}{2}}\Big(\mc D^{[k]}({\bm t})\mc D^{[n-k]}({\bm t})  \Big)^{m_k(\pi)/2}.
 \eeq
 \item For  $q \equiv 3 \pmod 4$,
 \beq
 w(\pi)=\prod\limits_{k=1}^{\frac{n-1}{2}}\Big(\mc D^{[k]}({\bm v}) \mc D^{[n-k]}({\bm v}) \Big)^{m_k(\pi)}. 
 \eeq
 \end{roster}
\end{thm}

Here as in Section \ref{revGJ}, $\bm t=\{t_1, \ldots, t_{n-1}\}$ with $t_i \in H^*(T)$, and $\mc D^{[k]}=\mc D_{n-1}^{[k]}$ when $q \equiv 1 \pmod 4$, and similarly for $q \equiv 3 \pmod 4$.
Note that $m_k(\pi)$ is given by Lemma \ref{mkpi}.

\begin{proof} The restriction of $\pi$ to $\GL(n-1,q)$ has $\det \pi=1$. The theorem thus follows from formulas \eqref{q34}, \eqref{q14} together with \eqref{office}. 
\end{proof} 

\begin{remark} See \cite{Malik.thesis} for a direct proof, which does not make use of \cite{GJgln}. (In fact, both these papers grew out of the same seminar.)
\end{remark}

\begin{ex}
 If $\pi$ is the regular representation of $G$, then $\chi_\pi(b_i)=0$ for $i>0$, so 
\beq
m_k(\pi)=m:=\frac{1}{2^{n-1}}|G|
\eeq
for all $k$. Therefore
\beq
w(\pi)= \begin{cases}   \mc D(\bm t)^{m/2} &\text{ when $q \equiv 1 \pmod 4$} \\
	            \mc D(\bm v)^{m} &\text{ when $q \equiv 3 \pmod 4$}. \\
	\end{cases}
	\eeq
\end{ex}

\section{The Case of $\SL(3,q)$}\label{sl33}
Let us now specialize Theorem \ref{SLnSWCs} to $G=\SL(3,q)$. 
First note that Lemma \ref{mkpi} gives $m_1(\pi)=\frac{1}{4}(\deg \pi-\chi_\pi(b_2))$, where
 
\beq
b_2=
   \begin{pmatrix}
      -1 & & \\
       & -1 & \\
       &&1\\
    \end{pmatrix}.
\eeq
For two variables, we have $\mc D =\mc D^{[1]} \mc D^{[2]} $, so that
\beq
w(\pi) =\begin{cases} \mc D(\bm t)^{\frac{m_1(\pi)}{2}}, & \text{ for $q \equiv 1 \pmod{4}$} \\
\mc D(\bm v)^{m_1(\pi)}, & \text{ for $q \equiv 3 \pmod{4}$}. \\
	\end{cases}
	\eeq

\begin{lemma} \label{good.friday} The integer $m_1(\pi)$ is even.
\end{lemma}
\begin{proof}
When $q \equiv 1 \pmod 4$, this follows from Lemma \ref{lien}.  For $q \equiv 3 \pmod 4$,   the above gives
$w_2(\pi)=m_1(\pi)(v_1^2+v_2^2+v_1v_2)$, which vanishes by Lemma \ref{Schur}.
 \end{proof}

So write $ e(\pi)=\half m_1(\pi)=\cfrac18\big(\deg \pi-\chi_\pi(b_2)\big)$. 
For $i=1,2$, put $$d_i=\begin{cases}
d_i(t_1,t_2) & q\equiv  1\: (\Mod\: 4)\\
d_i(v_1^2,v_2^2) & q\equiv  3\: (\Mod\: 4).
\end{cases}$$
Note that $\deg d_1=4$ and $\deg d_2=6$. Also put $\mc D=1+d_1+d_2$. From Theorem \ref{SLnSWCs}, we simply have:
\begin{equation} \label{Saturday}
w(\pi) =\mc D^{e(\pi)} =(1+d_1+d_2)^{e(\pi)}. 
\end{equation}

In particular, $w_4(\pi)=e(\pi) d_1$ and $w_6(\pi)=e(\pi)d_2$.
 
\begin{cor}\label{sl3obs}
  The obstruction class of an orthogonal representation $\pi$ of $G$ is $d_1^{2^r}$, where
  $r=\ord_2(e(\pi))$.  
\end{cor}

\begin{proof}
From \eqref{Saturday} we have
$$w(\pi)=\sum_{i=0}^{e(\pi)}{e(\pi) \choose i}(d_1+d_2)^i.$$
By elementary number theory, $e(\pi) \choose2^r$ is the first odd binomial coefficient appearing in the above sum.  Therefore the obstruction class is $d_1^{2^r}$, as claimed.
\end{proof}

\begin{cor}\label{HSW} We have $H^*_{\SW}(G)=\z/2\z[d_{1},d_{2}]$. Moreover, the Stiefel-Whitney group $\mc W(G)$ is infinite cyclic with generator $\mc D$.
\end{cor}

\begin{proof}
We will produce two orthogonal representations $\Pi_1$ and $\Pi_2$ with $e(\Pi_1)$ and $e(\Pi_2)$ coprime. This will suffice to prove the second assertion. One of these two integers must be odd, say $e(\Pi_1)$. 
Then $w_4(\Pi_1)=d_1$ and $w_4(\Pi_1)=d_2$, and this gives the first assertion.
  
Let $B$ be the subgroup of upper triangular matrices, and $L$ be the following Levi subgroup of $G$:
$$L=\Bigg\{\begin{pmatrix}
\ast&0&0\\
0&\ast&\ast\\
0&\ast&\ast
\end{pmatrix}\Bigg\}.$$
Note that $L$ is isomorphic to $\GL(2,q)$.
Let $P$ be the parabolic subgroup generated by $L$ and $B$. 

Let $\chi$ be a linear character of  $\f_q^\times$ with $\chi(-1)=-1$, and let $\alpha$ be the linear character $\chi \circ \det$ of $L$. Then $\alpha(\diag(-1,-1,1))=-1$.
We  consider the parabolic induction $\pi_1=\Ind_P^G\as$. Put $\Pi_1=S(\pi_1)$. We have
$\deg(\pi_1)=q^2+q+1$ and $\chi_{\pi_1}(b_2)=-q$, whence $e(\Pi_1)=\frac14(q^2+2q+1)$.
[Here and below we use the Frobenius character formula for induced representations and identified the $G$-set $G/P$ with the projective plane $\mb P^2(\mb F_q)$. Note that $b_2$ acts on $\mb P^2(\mb F_q)$ by $(x:y:z)\mapsto (-x:-y:z)$.]

Let $\varphi$ a linear character of $\f_{q^2}^\times$ with $\varphi(-1)=-1$, and write $\sigma_\varphi$ be the corresponding cuspidal representation of $L$. (In \cite[Section 5.2]{Fulton.Harris} this is called $X_\varphi$.) Such a representation has degree $(q-1)$. Let $\pi_2$ be the parabolically induced representation $\pi_2=\Ind_P^G\sigma_\varphi$, and put $\Pi_2=S(\pi_2)$. Now $\deg\pi_2=(q-1)(q^2+q+1)$ and $\chi_{\pi_2}(b_2)=1-q$, whence $e(\Pi_2)=\frac14(q-1)(q^2+q+2)$.    
 
The integers $\frac14(q^2+2q+1)$ and  $\frac14(q-1)(q^2+q+2)$ are clearly coprime.

\end{proof}

 Now Corollary \ref{G.topSWC} and \eqref{m0.simp} gives:

\begin{cor} \label{eggs} The class $w_{\Top}(\pi)$ vanishes unless $\deg\pi=-3\chi_\pi(b_2)$, in which case $w_{\Top}(\pi)=d_2^{e(\pi)}$. 
\end{cor}
 
 This curious condition invites further study. We offer the following consequence:
 
 \begin{cor} If $q \equiv 2 \mod 3$, then $w_{\Top}(\pi)=0$ for any orthogonal representation $\pi$ of $\SL(3,q)$.
 \end{cor}
 
 \begin{proof}  Say $\pi$ is irreducible. By Exercise 6.9 in \cite{serre.linear}, the quantity
   \beq
        \frac{\chi_\pi(b_2)}{\deg \pi}\mc C(b_2)
    \eeq
    is an integer. Here $\mc C(b_2)$ is the size of the conjugacy class of $b_2$. Since the centralizer of $b_2$ is $M \cong \GL(2,q)$, we have $\mc C(b_2)=q^2(q^2+q+1) \equiv 1 \mod 3$.
Now by   Corollary \ref{eggs},  if $w_{\Top}(\pi) \neq 0$, then $\mc C(b_2)$ would be a multiple of $3$, a contradiction.  It is a routine matter to extend this to the case of $\pi$ reducible.
  \end{proof}

 \section{The Case of $\SL(5,q)$}\label{sl5}

  Let $G=\SL(5,q)$ with  $b_2=\diag(-1,-1,1,1,1),$ $b_4=\diag(-1,-1,-1,-1,1)\in G$. Let $\pi$ be an orthogonal representation of $G$. 
  \subsection{Total SWC}
  
  From Theorem \ref{SLnSWCs}, the total SWC of $\pi$ is
\begin{equation}\label{sl5wpi}
w(\pi)=\begin{cases}
 ( \mc D^{[1]}(\bm t) \mc D^{[4]}(\bm t))^{m_1/2}  ( \mc D^{[2]}(\bm t) \mc D^{[3]}(\bm t))^{m_2/2}  & q \equiv 1 \pmod 4\\
( \mc D^{[1]}(\bm v) \mc D^{[4]}(\bm v))^{m_1}  ( \mc D^{[2]}(\bm v) \mc D^{[3]}(\bm v))^{m_2}  & q \equiv 3 \pmod 4,

\end{cases}
\end{equation}
 where
\begin{align*}
m_1&=\frac{1}{16}\big(\chi_\pi(\mathbbm{1})+2\chi_\pi(b_2)-3\chi_\pi(b_4)\big) \text{ and }\\
m_2&=\frac{1}{16}\big(\chi_\pi(\mathbbm{1})-2\chi_\pi(b_2)+\chi_\pi(b_4)\big).
\end{align*}
From Section \ref{dsonfac}, we have
 \beq
 \mc D^{[1]} \mc D^{[4]}=1+ (\mc E_1^2+\mc E_2)+ (\mc E_1 \mc E_2+\mc E_3)+ (\mc E_4+ \mc E_1 \mc E_3)+ \mc E_1 \mc E_4.
 \eeq
 and
 \begin{equation*}
\begin{split}
 \mc D^{[2]} \mc D^{[3]}&= 1+ (\mc E_1^2+\mc E_2)+(\mc E_1\mc E_2+\mc E_3)+(\mc E_1^4+\mc E_2^2 + \mc E_1 \mc E_3+\mc E_4)+\mc E_1 \mc E_4\\
&\quad+ (\mc E_1^6+\mc E_1^4\mc E_2+\mc E_2^3+\mc E_3^2)+(\mc E_1^5\mc E_2+\mc E_1\mc E_2^3+\mc E_1^4\mc E_3+\mc E_2^2\mc E_3) \\
&\quad+ (\mc E_1^4\mc E_2^2+\mc E_1^2\mc E_2^3+\mc E_1^5\mc E_3+\mc E_1\mc E_2^2\mc E_3+\mc E_1^2\mc E_3^2+\mc E_2\mc E_3^2+\mc E_1^4\mc E_4+\mc E_1^2\mc E_2\mc E_4+\mc E_2^2\mc E_4\\
&\quad\quad+\mc E_1\mc E_3\mc E_4) + \cdots,
 \end{split}
\end{equation*}
 where we have omitted the terms of degrees $9$ and $10$.
 \begin{remark}
By Corollary \ref{G.topSWC}, $w_{\Top}(\pi)$ is non-zero iff $\deg\pi=-10\chi_\pi(b_2)-5\chi_\pi(b_4)$.
\end{remark}
  \subsection{First few SWCs} \label{SL5Q3}
For nonnegative integers $A_1,A_2$, put
\beq
f=f^{A_1,A_2}= ( \mc D^{[1]}  \mc D^{[4]})^{A_1}  ( \mc D^{[2]} \mc D^{[3]} )^{A_2} \in \z/2\z[\mc E_1, \ldots, \mc E_4].
\eeq
Write $f_k$ for the degree $k$ term of $f$.
We have $f_1=0$ and $f_2=(A_1+A_2)\mc H_2$. (Recall that $\mc H_2=\mc E_1^2+\mc E_2$.) Hence $f_2=0$ iff $A_1 \equiv A_2 \pmod 2$.
When this holds, then
\beq
f_4=\left(\binom{A_1}{2}+\binom{A_2}{2}\right)\mc H_2^2,
 \eeq
and this vanishes iff $A_1 \equiv A_2 \pmod 4$. When the latter holds, then (by a long calculation):
\beq
f_8=\left(\binom{A_1}{4}+\binom{A_2}{4} \right) \mc H_2 ^4+A_1 \left(\mc E_1 ^2\mc E_2 \mc E_4 +\mc E_1 \mc E_3 \mc E_4 +\mc E_4 ^2+\mc E_1 ^2\mc E_3 ^2+\mc H_2^4 \right).
\eeq

We have
\beq
w(\pi)=\begin{cases} f^{\half m_1,\half m_2}(\bm t) &\text{ when $q \equiv 1 \pmod 4$}, \\
				f^{m_1,m_2}(\bm v)  &\text{ when $q \equiv 3 \pmod 4$}. \\
				\end{cases}
				\eeq
  From this we immediately obtain formulas for $w_4(\pi)$, $w_8(\pi)$, and $w_{16}(\pi)$ (subject to vanishing conditions) when  $q \equiv 1 \pmod 4$, and $w_4(\pi)$ and $w_8(\pi)$ when $q \equiv 3 \pmod 4$. For example, when $q \equiv 1 \pmod 4$ and $m_1 \equiv m_2 \pmod 4$, then
  $w_4(\pi)=0$ and 
\beq
w_8(\pi)=\left(\binom{\half m_1}{2}+\binom{\half m_2}{2}\right) \mc H_2(\bm t)^2.
\eeq
  
 When $q \equiv 1 \pmod 4$, supposing $m_1 \equiv m_2 \pmod 8$, then
 \beq
 w_{16}(\pi)=\left[\left(\binom{\half m_1}{4}+\binom{\half m_2}{4} \right) \mc H_2 ^4+\half m_1 \cdot \left(\mc E_1 ^2\mc E_2 \mc E_4 +\mc E_1 \mc E_3 \mc E_4 +\mc E_4 ^2+\mc E_1 ^2\mc E_3 ^2+\mc H_2^4 \right)\right](\bm t).
 \eeq
 When $q \equiv 3 \pmod 4$ and $m_1 \equiv m_2 \pmod 4$, then
 \beq
 w_{8}(\pi)=\left[\left(\binom{m_1}{4}+\binom{m_2}{4} \right) \mc H_2 ^4+ m_1 \cdot \left(\mc E_1 ^2\mc E_2 \mc E_4 +\mc E_1 \mc E_3 \mc E_4 +\mc E_4 ^2+\mc E_1 ^2\mc E_3 ^2+\mc H_2^4 \right) \right](\bm v).
 \eeq

\section{Universal SWCs}\label{uniswc}

 Let $G_{n}=\SL(n,q)$. For $n_0 < n$, we may view $G_{n_0}$ as subgroup of $G_n$ via $ i_{n_0}(g)=\begin{pmatrix} g & \\  &I \\ \end{pmatrix}$. In this section we give a detection theorem for $G_{n_0} \leq G_{n}$, which allows for a ``universal'' calculation
  of $w_k(\pi)$ in terms of character values, valid for all $n$ sufficiently large.

 \subsection{Detection by Smaller Special Linear Groups}
  
\begin{prop}\label{Nakao}
 Let $n,n_0$ be odd with $n \geq n_0$. The restriction map
$i^*_{n_0}:H^k(G_{n})\to H^k(G_{n_0}) $
is injective for all $k<n_0$ when $q\equiv 3\;(\Mod\; 4)$ and for $k<2n_0-1$ when $q\equiv 1\;(\Mod\; 4)$.
\end{prop}

\begin{proof}
As before, we have   subgroups
$$M_n=\Big\{\begin{pmatrix}
A&&0\\
0&&\det(A)^{-1}
\end{pmatrix}:A\in \GL(n-1,q)\Big\}.$$
The diagram
\beq
	\xymatrix{
		M_{n_0} \ar[r]  \ar[d] & G_{n_0} \ar[d]^{i_{n_0}}\\
		M_n\ar[r] & G_n\\
	}
	\eeq 
does \emph{not} commute: If $g \in M_{n_0}$, then going right, then down gives
\begin{equation} \label{green}
\begin{pmatrix} g && \\
			& (\det g)^{-1} & \\
			&& I,
			\end{pmatrix},
\end{equation}
but going down, then right gives
\begin{equation} \label{black}
\begin{pmatrix} g && \\
			& I & \\
			& &(\det g)^{-1} \\
			\end{pmatrix}.
\end{equation}
However there is an evident permutation matrix $A$ (of order $2$) conjugating \eqref{green} to \eqref{black}. Write $\theta$ for the automorphism of $G_n$ induced by $A$. Then, the following diagram commutes:
 
 $$\begin{tikzcd} 
		&H^k(G_{n}) \ar[r,"i^*_{n_0}"]   & H^k(G_{n_0}) \ar[d,hook]\\
		H^k(G_n)\ar[r,hook] \ar[ur,"\theta^*"]& H^k(M_n)\ar[r,hook]& H^k(M_{n_0})&\\
 \end{tikzcd}$$

The first map in the lower row is injective by Lemma \ref{dl}, as is the second map by \cite[Theorems 3 and 8]{GJgln} all $k<n_0$ when $q\equiv 3\;(\Mod\; 4)$ and for $k<2n_0-1$ when $q\equiv 1\;(\Mod\; 4)$. Since $\theta^*$ is an automorphism, $i_{n_0}^*$ is injective, as claimed.
\end{proof}

\begin{remark} \label{univ.remark} When $q \equiv 3 \pmod 4$, and $k$ is even, the SWC $w_k(\pi)$ can be computed by restricting to $\SL(k+1,q)$. When $q \equiv 1 \pmod 4$ and $k$ is a multiple of $4$, then we may compute $w_k(\pi)$ by restricting to $\SL \left(\frac{k}{2}+1,q\right)$.
\end{remark}

\subsection{Universal $w_4$ and $w_8$} 
 \begin{thm}\label{w4calc} Let $\pi$ be an orthogonal representation of $\SL(n,q)$ with $n \geq 5$ odd.  
Put  \begin{align*}
m_1&=\frac{1}{16}\big(\chi_\pi(\mathbbm{1})+2\chi_\pi(b_2)-3\chi_\pi(b_4)\big) \text{ and }\\
m_2&=\frac{1}{16}\big(\chi_\pi(\mathbbm{1})-2\chi_\pi(b_2)+\chi_\pi(b_4)\big).
\end{align*}

\begin{roster}
\item When $q \equiv 3 \pmod 4$,
\beq
w_4(\pi)=\left( \binom{m_1}{2}+\binom{m_2}{2} \right) \mc H_2(\bm v)^2.
\eeq
 
 \item 
When $q \equiv 1 \pmod 4$, we have $w_4(\pi) =\half(m_1+m_2) \mc H_2(\bm t)$. When also $w_4(\pi)=0$, 
 \beq
w_8(\pi)=\left( \binom{\half m_1}{2}+ \binom{\half m_2}{2} \right) \mc H_2(\bm t)^2.
\eeq

\end{roster}
\end{thm}

\begin{proof}
\begin{roster}\item When  $q\equiv 3\pmod 4$, the restriction $i^*_5:H^4(G_n)\to H^4(G_5).$ is injective by Proposition \ref{Nakao}. So the conclusion for $w_4(\pi)$ follows from Section \ref{SL5Q3}.

\item When $q \equiv 1 \pmod 4$,  the restriction maps from $H^4(G_n)$ and $H^8(G_n)$ to $H^4(G_5)$ and $H^8(G_5)$ are injective,  so we may use the formulas in Section \ref{SL5Q3}.
  \end{roster}
\end{proof}

Clearly Theorem \ref{intro.thm} follows from this.

\bibliographystyle{alpha}
\bibliography{mybib}
\vspace{10mm}

\end{document}